\documentclass[a4paper,10pt,reqno]{amsart}

\usepackage{graphicx}
\usepackage{a4wide}
\usepackage{amssymb}
\usepackage{amsthm}
\usepackage{eucal}
\usepackage[pdftex,colorlinks]{hyperref}

\numberwithin{equation}{section}

\newcommand{\R}{\mathbb R}

\newcommand{\be}{\begin{equation}}
\newcommand{\ee}{\end{equation}}
\newcommand{\ba}{\begin{eqnarray}}
\newcommand{\ea}{\end{eqnarray}}

\newcommand{\pt}{\ensuremath{\partial_t}}
\newcommand{\norm}[1]{\ensuremath{\left\|#1\right\|}}
\newcommand{\Div}{\mathrm{div}\,}
\newcommand{\bM}{\ensuremath{\mathbf{M}}}
\newcommand{\Grad}{\mathrm{\nabla}}
\newcommand{\eps}{\ensuremath{\varepsilon}}
\newcommand{\abs}[1]{\ensuremath{\left|#1\right|}}
\newcommand{\pxij}{\ensuremath{\partial_{x_ix_j}^2}}
\newcommand{\pxi}{\ensuremath{\partial_{x_i}}}
\newcommand{\pxj}{\ensuremath{\partial_{x_j}}}
\newcommand{\J}{\ensuremath{\mathcal{J}}}
\newcommand{\dx}{\ensuremath{\, dx}}
\newcommand{\dt}{\ensuremath{\, dt}}

\newtheorem{theorem}{Theorem}[section]

\newtheorem{remark}[theorem]{Remark}
\newtheorem{question}[theorem]{Question}
\newtheorem{lemma}[theorem]{Lemma}

\begin{document}

\title[Uniform Controllability for a Degenerating System]{Uniform Null Controllability for a Degenerating Reaction-Diffusion System Approximating a Simplified Cardiac Model}

\author[F. W. Chaves-Silva]{Felipe Wallison Chaves-Silva$^*$}
\address[F. W. Chaves-Silva]{\newline Universit\'e de Nice Sophia-Antipolis, Laboratoire Jean Dieudonn\'e, UMR CNRS 6621, Parc Valrose, 06108 Nice Cedex 02,
France.}
\email{fchaves@unice.fr}\thanks{$^*$F. W. Chaves-Silva has been supported  by  the ERC project Semi Classical Analysis of  Partial Differential Equations, ERC-2012-ADG, project number 320845;  the Grant  BFI-2011-424 of the Basque Government and partially supported by the Grant  MTM2011-29306-C02-00 of the MICINN, Spain, the ERC Advanced Grant FP7-246775
NUMERIWAVES, ESF Research Networking Programme OPTPDE and the Grant PI2010-04 of the Basque Government}

\author[M. Bendahmane]{Mostafa Bendahmane}
\address[M. Bendahmane]{\newline
Institut Math\'ematiques de Bordeaux, Universit\'e Victor Segalen Bordeaux 2, 3 ter Place de la Victoire, 33076 Bordeaux,
France.}
\email[]{mostafa.bendahmane@u-bordeaux2.fr}



\begin{abstract} 
This paper is devoted to the analysis of the uniform null controllability  for a family of nonlinear reaction-diffusion systems
approximating a parabolic-elliptic system which models the  electrical activity of the heart. The uniform,
with respect to the degenerating parameter,  null controllability of the approximating system by means of a
single control is shown. The proof is based on the combination of Carleman estimates and weighted energy inequalities.
\end{abstract}

\maketitle

\section{Introduction}
\label{sec:intro}

Let $\Omega \subset \mathbb{R}^N$ $(N= 2, 3)$ be a bounded connected open set whose boundary, $\partial \Omega$, is sufficiently regular. Let $T >0$, and let $\omega$ and $\mathcal{O}$ be  two (small) nonempty subsets of $\Omega$, which we will refer  to as \textit{control domains}. We will use the notation $Q = \Omega \times (0,T)$ and $\Sigma = \partial \Omega \times (0,T)$.

The main objective of this paper is to study the  properties of controllability and observability for a  family of nonlinear reaction-diffusion systems which degenerates into a nonlinear parabolic-elliptic system which models the electrical activity in the cardiac tissue.

To state the model, we let $u_i=u_i(t,x)$ and  $u_e=u_e(t,x)$ represent the \textit{intracellular} and
\textit{extracellular} electric potentials, respectively. Their difference,
$v=u_i-u_e$, is called the \textit{transmembrane} potential. The
anisotropic properties of the media are modeled by intracellular and
extracellular conductivity tensors $\bM_i(x)$ and $\bM_e(x)$.

The widely accepted model (see \cite{Colli4, Henr,Tung}) describing the electrical activity in the cardiac tissue reads as follows:
\begin{equation}\label{S1}
\left\{
\begin{array}{llcc}
c_m \pt v -\Div (\bM_i(x)\Grad u_i) +h(v)= f1_{\omega} &\qquad \text{ in } Q,  \\ 
 c_m \pt v  + \Div( \bM_e(x)\Grad u_e ) +h(v )=g1_{\mathcal{O}} &\qquad \text{ in } Q,
 \end{array}
 \right.
\end{equation}
where $c_m>0$ is the  surface capacitance of the membrane, the nonlinear function $h: \mathbb{R} \rightarrow \mathbb{R}$ is the transmembrane ionic current (the most interesting case being when $h$ is a cubic polynomial), and $f$ and $g$ are stimulation currents applied, respectively, to $\omega$ and $\mathcal{O}$. 
\\
System (\ref{S1}) is  known as the  \textit{bidomain model} and is completed with Dirichlet boundary 
conditions for the intra- and extracellular electric potentials
\begin{equation}
	\label{S2}
	u_i= u_e=0 \  \text{on} \ \Sigma \\
\end{equation}
and initial data for the transmembrane potential
\begin{equation}
	\label{S3}
	v(0,x)=v_0(x), \quad x\in \Omega.
\end{equation}
We point out that realistic models describing electrical activities in the heart also include a system of ODE's for computing the ionic current as a function of the transmembrane potential and a series of additional ``gating variables'' which aim  to model the ionic transfer across the cell membrane (see  \cite{Hodgkin-Huxley,Keener, LuoRudy,Noble}).

In the case where $f1_\omega = g1_{\mathcal{O}}$ and   $\bM_i = \mu \bM_e$, for some constant $\mu \in \R$, the bidomain model is simplified into the following parabolic-elliptic system:
\begin{equation}
	\label{Sequiv1}
	\begin{cases}	
 c_{m} \pt v-\frac{\mu}{\mu+1}\Div\bigl(\bM_{e}(x)\Grad v \bigr)+
h(v) = f1_{\omega}  &\qquad \text{ in } Q, \\
 -\Div\bigl(\bM(x)\Grad u_{e}\bigr)  = \Div\bigl(\bM_i(x)\Grad v\bigr) &\qquad \text{ in } Q, \\
 v = u_e = 0  & \qquad \text{ on } \Sigma, \\
  v(0) = v_0  &\qquad  \text{ in } \ \Omega,
	\end{cases}	
\end{equation}
where $M = M_i + M_e$. 

System  \eqref{Sequiv1} is known as \textit{monodomain model} and is a very interesting model from the implementation point of view, since it conserves some of the essential features of the bidomain model as excitability phenomena  (see \cite{Colli4, Kunisch, Veneroni}). 

The main difference between the  bidomain model \eqref{S1} and the monodomain model \eqref{Sequiv1} is the fact that the first model is a system of two coupled parabolic equations, while the second one is a system of parabolic-elliptic type. Therefore, from the control point of view, one could expect these two systems to have, at least a priori, different control properties. In this work  we show that the properties of controllability and observability for the monodomain model can be seen as  a limit process of  the  controllability properties  of a  family of  coupled parabolic systems. Indeed, given $\varepsilon \in \mathbb{R}$ such that $0< \varepsilon \leq 1$, we approximate the monodomain model by the following  family of parabolic systems:
\begin{equation}
	\label{Srelaxed-orig}
	\begin{cases}	
 c_{m} \pt v^{\eps}-\frac{\mu}{\mu+1}\Div\bigl(\bM_{e}(x)\Grad v^{\eps} \bigr)+
h(v^{\eps}) = f^{\eps}1_{\omega} &\qquad \text{ in } Q, \\
 \eps \pt u_e^{\eps}-\Div\bigl(\bM(x)\Grad u^{\eps}_{e}\bigr)  = \Div\bigl(\bM_{i}(x)\Grad v^{\eps}\bigr) 
 &\qquad \text{ in } Q, \\
 v^{\eps} = u^{\eps}_e = 0  &\qquad \text{ on } \Sigma, \\
  v^{\eps}(0) = v_0, \ u_e^{\eps}(0) = u_{e,0} &\qquad  \text{ in } \ \Omega.
	\end{cases}	
\end{equation}
In this paper we give a  positive  answer to the following question:

\begin{question}\label{q1} If, for each $\varepsilon >0$,  there exists a control $f^{\varepsilon}$ that drives  the solution $(v^{\varepsilon}, u_e^{\varepsilon})$ of (\ref{Srelaxed-orig}) to zero at time $t= T$, i.e.,
$$
v^{\varepsilon}(T) = u_e^{\varepsilon}(T) = 0,
$$
is it true that,  when $\varepsilon \rightarrow 0^+$, the control sequence $\{f^{\varepsilon}\}_{\varepsilon >0}$ converges to a function $f$ which drives  the associated solution  $(v,u_e)$ of (\ref{Sequiv1}) to zero at time $t=T$?
\end{question}

This question of approximating an equation by another having different physical properties has been used several times in the case of parabolic equations degenerating into hyperbolic ones (see, for example,  \cite{Coron-1, Glass, Gue-2}) and hyperbolic equations degenerating into parabolic ones (see, for example,  \cite{Lopez-1, Lopez-2}). However, as far as we know, this is the first time that controllability of parabolic systems degenerating into parabolic-elliptic systems
is studied. It is also important to mention that families of parabolic systems which degenerate
into parabolic-elliptic ones arise in many areas,  such as biology, chemistry and astrophysics (see  \cite{Biler-2,Biler-3, Keller-Segel}).

As usual, in control theory, when dealing with the controllability of a nonlinear problem,  we  first consider the linearized version of  (\ref{Srelaxed-orig}):
\begin{equation}
	\label{Srelaxed-orig-1}
	\begin{cases}
c_{m} \pt v^{\varepsilon}-\frac{\mu}{\mu+1}\Div\bigl(\bM_{e}(x)\Grad v^{\varepsilon} \bigr)+
a(t,x)v^{\varepsilon} = f^{\varepsilon}1_{\omega} &\qquad \text{ in } Q, \\
\eps \pt u^{\varepsilon}_e-\Div\bigl(\bM(x)\Grad u^{\varepsilon}_{e}\bigr)  = \Div\bigl(\bM_{i}(x)\Grad v^{\varepsilon}\bigr)&\qquad \text{ in } Q, \\
 v^{\varepsilon} = u^{\varepsilon}_e = 0  &\qquad \text{ on } \Sigma, \\
 v^{\varepsilon}(0) = v_0, \ u^{\varepsilon}_e(0) = u_{e,0} &\qquad  \text{ in } \ \Omega,
	\end{cases}	
\end{equation}
where $a$ is a bounded function.

Given $\varepsilon >0$, the first obstacle to answering, positively, Question \ref{q1}, will be to drive  $(v^{\varepsilon},u^{\varepsilon}_e)$, solution of (\ref{Srelaxed-orig-1}), to zero at time $T$ by means of a control $f^{\varepsilon}$ in such a way  that the  sequence of controls $\{f^{\varepsilon}\}_{\varepsilon >0}$ converges when $\varepsilon \rightarrow 0^+$. Once it is shown that  such a convergent sequence of control, $\{f^{\varepsilon}\}_{\varepsilon >0}$, for the linear system \eqref{Srelaxed-orig-1}, exists, we employ a fixed point argument and conclude that the same is true for the nonlinear system  (\ref{Srelaxed-orig}). 

Thus, we introduce the adjoint system of (\ref{Srelaxed-orig-1}):
 \begin{equation}
	\label{generaladjoint-1}
	\begin{cases}
		 - c_m \pt \varphi^{\eps} -\frac{\mu}{\mu +1}\Div (M_e(x)\Grad \varphi^{\eps}) 
		+a(t,x)\varphi^{\eps}= \Div (M_i(x)\Grad \varphi_e^{\eps})& \text{ in } Q ,\\ 
		 - \eps \pt \varphi^{\eps}_{e}- \Div (M(x)\Grad \varphi^{\eps}_e)=0& \text{ in } Q, \\
		\varphi^{\eps}= \varphi^{\eps}_e = 0 & \text{ on } \Sigma, \\
		 \varphi^{\eps}(T) = \varphi_T, \  \varphi^{\eps}_e(T) = \varphi_{e,T} & \text{ in } \ \Omega.
	\end{cases}	
\end{equation}

Using duality arguments, it is very easy to prove that the task of building such a convergent sequence of controls, $\{f^{\varepsilon}\}_{\varepsilon >0}$, for \eqref{Srelaxed-orig-1}  is  equivalent to prove the  following (uniform) observability inequality for the solutions of \eqref{generaladjoint-1}:
\begin{equation}\label{observability-inequality}
||\varphi^{\eps}(0)||^2_{L^2(\Omega)} + \varepsilon ||\varphi^{\eps}_e(0)||_{L^2(\Omega)}^2  \leq C\iint_{Q_{\omega}}|\varphi^{\eps}|^2  dxdt,
\qquad Q_{\omega}:= \omega \times (0,T),
\end{equation}
where $(\varphi_{T}, \varphi_{e,T}) \in L^2(\Omega)^2$ and the constant $C = C(\varepsilon, \Omega, \omega, ||a||_{L^\infty},  T)$ remains bounded when $\eps \rightarrow 0^+$. 

We prove inequality  \eqref{observability-inequality} as a consequence of an appropriate Carleman inequality for the solution $(\varphi^{\eps}, \varphi^{\eps}_e)$  of \eqref{generaladjoint-1}  (see section \ref{CI}). 
We notice that, due to the fact the control is acting on the first  equation of \eqref{Srelaxed-orig-1}, in our Carleman inequality, we need to bound global integrals of $\varphi^{\eps}$ and $\varphi^{\eps}_{e}$ in terms of a local integral of $\varphi^{\eps}$, uniformly with respect to $\varepsilon$. Two main difficulties appear: first, the coupling in the first equation of \eqref{generaladjoint-1} is in $\Div (M_i(x)\Grad \varphi_e^{\eps})$ and not in $\varphi_e^{\eps}$; second, we must  show that the constant we get in our  Carleman inequality does not blow up when $\varepsilon \rightarrow 0^+$.

The first difficulty is not so hard to overcome. Indeed,  for each $\varepsilon >0$ fixed, inequality \eqref{observability-inequality} is known to be true for system \eqref{generaladjoint-1} (see \cite{Gue-1}). However, the main novelty here is the fact that we obtain the boundedness of the observability constant $C$ with respect to $\varepsilon$. As we will see,  Carleman inequalities alone are not enough for this task, and we need to combine sharp Carleman estimates, with respect to $\varepsilon$, and weighted energy inequalities.
 
As far as the controllability of non degenerate coupled parabolic systems is concerned, the situation is, by now, fairly well understood. For instance, in \cite{Gue-1}, the controllability of a quite general linear coupled parabolic system is studied and a null controllability result is obtained by means of  Carleman inequalities.   In \cite{Khodja-1}, using a different strategy,   the controllability of a
nonlinear reaction-diffusion system of  two coupled parabolic equations is analyzed, and the authors prove the null controllability for the linear system  and the local null controllability of the nonlinear one.  Another relevant work concerning the controllability of coupled systems is  \cite{Burgoz-deTeresa}, in which the authors analyze
the null controllability of a cascade system of $m$ ($m > 1$) coupled parabolic equations and the authors
are able to obtain null controllability for the cascade system, whenever they have a  good coupling structure.  It is also worth mentioning the works  \cite{Khodja-2}, \cite{Barbu} and  \cite{Burgos-Garcia}, where local and global controllability results for phase field systems were studied. 

For a general discussion about the controllability of coupled parabolic systems, see the survey paper \cite{Kho-Be-Bu-Te}.

Concerning controllability results for the  bidomain model, since in both equations the couplings are given by the time derivatives of the electrical potentials, it seems very difficult to study controllability properties for such a model.  To the  best of our  knowledge, for the bidomain model \eqref{S1},  the problems of null and approximate controllability   are still open (even with two controls). Regarding the null controllability of the  monodomain model (\ref{Sequiv1}), since the solution of the parabolic equation enters as a source term in the elliptic one, the following controllability result holds.

\begin{theorem}\label{controllability3}
 \begin{enumerate}
 \item  If $ h$ is  $C^1(\mathbb{R})$,  globally Lipschitz and $h(0)= 0$, then, for every  $v_0 \in L^2(\Omega)$, there exists a control $f \in L^2(\omega \times (0,T))$ such that the solution $(v,u_e)$ of  (\ref{Sequiv1}) satisfies:
$$
v(T) = u_e(T) = 0.
$$
Moreover, the control $f$ satisfies the following estimate:
\begin{equation}\label{bound-control-1}
\norm{f1_{\omega}}^2_{L^2(Q)} \leq C \norm{v_0}^2_{L^2(\Omega)},
\end{equation}
for  a constant $C = C(\Omega, \omega,  T)>0$.

\item  If $h$ is $C^1(\mathbb{R})$ and  $h(0)= 0$,  there exists $\gamma >0$ such that, for every  $v_0 \in W^{2/3, 6}(\Omega)\cap H^1_0(\Omega)$ with  $||v_0||_{W^{2/3, 6}(\Omega)\cap H^1_0(\Omega)} \leq \gamma $, there exists a control $f \in L^{6}(\omega \times (0,T))$ such that the solution $(v,u_e)$ of  (\ref{Sequiv1}) satisfies:
$$
v(T) = u_e(T) = 0.
$$
Moreover, the control $f$ satisfies the following estimate:
\begin{equation}\label{bound-control-1}
\norm{f1_{\omega}}_{L^{6}(Q)} \leq C \norm{v_0}_{L^2(\Omega)},
\end{equation}
for  a constant $C = C(\Omega, \omega,  T)>0$.
\end{enumerate}
\end{theorem}

Theorem \ref{controllability3}, case $1$, follows from \cite[Theorem 3.1]{16}  and case $2$ follows from \cite[Theorem 3.5]{Guerrero}   (see also \cite[Theorem 4.2]{16}).


This paper is organized as follows. In section \ref{main-result}, we state our main results.
In section \ref{CI}, we  prove a uniform Carleman inequality for the adjoint system (\ref{generaladjoint-1}). Next, we show, in  section \ref{control-linearsystem}, the uniform null controllability of (\ref{Srelaxed-orig-1}).
In section \ref{section-nonlinear},  we deal with the uniform null controllability of the nonlinear system (\ref{Srelaxed-orig}).



\section{Main results}
\label{main-result}
Throughout  this paper  we will assume that  the matrices $M_j$, $j=i,e$ are  $C^{\infty}$, 
bounded, symmetric and positive semidefinite. 

Our first main result is a uniform Carleman estimate for the adjoint system (\ref{generaladjoint-1}).

\begin{theorem}\label{thm-carleman-syst}
Given any $0< \epsilon \leq 1$, there exist positive constants $C = C(\Omega, \omega)$, $\lambda_0= \lambda_0(\Omega, \omega) \geq 1$ and $s_0= s_0(\Omega, \omega) \geq 1$ such that, for any  $(\varphi_T,  \varphi_{e,T}) \in L^2(\Omega)^2$ and any $  a \in L^{\infty}(Q)$, the solution  $(\varphi^{\eps}, \varphi_e^{\eps})$ of  (\ref{generaladjoint-1}) satisfies:
\begin{align}\label{observability1112}
 \int \! \! \! \int_{Q}e^{3s\alpha}| \rho^{\eps}|^2dxdt &+s^3\lambda^4\int \! \! \! \int_{Q}\phi^3e^{3s\alpha}|\varphi^{\eps}|^2dxdt \nonumber \\
 & \leq C e^{6\lambda||\psi||}s^8\lambda^4 \int \! \! \! \int_{Q_{\omega}}\phi^8 e^{2s\alpha}|\varphi^{\eps}|^2  dxdt, 
\end{align}
 for every $s\geq (T +(1+||a||_{L^{\infty}}^{2/3}) T^2)s_0$ and $\lambda \geq \lambda_0$, where  $\rho^{\eps}(x,t) = \Div (M(x)\nabla \varphi_e^{\eps}(x,t))$  and the
weight functions $\phi$ and $\alpha$ are defined in (\ref{weightfunctions1}) and (\ref{weightfunctions2}), respectively. 
\end{theorem} 

The proof of Theorem \ref{thm-carleman-syst} follows from a  combination of Carleman inequalities, for the heat equation, with a  precise dependence on the degenerating parameter, and an energy inequality for the adjoint system (\ref{generaladjoint-1}). We prove Theorem \ref{thm-carleman-syst}   in section \ref{CI}.

%
\begin{remark}
As a direct consequence of the Carleman inequality \eqref{observability1112},  we have the unique continuation
property for the solutions $(\varphi^{\eps}, \varphi^{\eps}_e)$ of   \eqref{generaladjoint-1}:
$$
\mbox{``Given} \ \varepsilon >0,  if  \ \varphi^{\eps} = 0 \ \mbox{in} \ \omega \times (0,T), \mbox{then} \ (\varphi^{\eps}, \varphi^{\eps}_e) \equiv (0,0) \  \mbox{in}  \ Q \text{''}.
$$
This unique continuation property for the adjoint system \eqref{generaladjoint-1} implies, for each $\varepsilon >0$, the approximate controllability at time $T$ of system \eqref{Srelaxed-orig-1},  with a control acting only on the first equation.
\end{remark}

The second main result of this paper gives the global null controllability of the linear system (\ref{Srelaxed-orig-1}).
\begin{theorem}\label{controllability2}
For any $0< \varepsilon \leq 1$ and any $(v_0, u_{e,0}) \in L^2(\Omega)^2$,  there exists a control $f^{\varepsilon} \in L^2(\omega \times (0,T))$ such that the associated solution, $(v^\eps, u^\eps_e)$, to   (\ref{Srelaxed-orig-1}) is driven to zero at time $T$. That is to say, the associated solution satisfies:
$$
v^\eps(T) = 0, \ u^\eps_e(T) = 0.
$$
Moreover, the control $f^{\eps}$ satisfies the estimate:
\begin{equation}\label{bound-control}
\norm{f^\eps1_{\omega}}^2_{L^2(Q)} \leq C\bigl( \norm{v_0}^2_{L^2(\Omega)} + \eps \norm{u_{e,0}}^2_{L^2(\Omega)}\bigl),
\end{equation}
for  a constant $C = C(\Omega, \omega, ||a||_{L^\infty},  T) > 0$.
\end{theorem}
\vspace{0.3cm}

From Theorem \ref{thm-carleman-syst}, the proof of Theorem \ref{controllability2} is standard. However, for the sake of completeness,  we prove Theorem \ref{controllability2}  in section \ref{control-linearsystem}.

The third main result of this paper is concerned with the uniform null controllability of the nonlinear parabolic system (\ref{Srelaxed-orig}).
\begin{theorem}\label{controllability4}
Given any $0< \varepsilon \leq 1$, we have: 
\begin{enumerate}
\item If $ h$ is  $C^1(\mathbb{R})$,   globally Lipschitz and $h(0)= 0$, then, for every $(v_0, u_{e,0}) \in L^2(\Omega)^2$, there exists a control $f^\eps \in L^2(\omega \times (0,T))$ such that the solution  $(v^{\eps},u^{\eps}_e)$ of   (\ref{Srelaxed-orig}) satisfies: 
$$
v^{\eps}(T) = u^{\eps}_e(T) = 0.
$$
Moreover, the control $f^{\eps}$ satisfies the estimate:
\begin{equation}\label{bound-control-1}
\norm{f^{\eps}1_{\omega}}^2_{L^2(Q)} \leq C\bigl( \norm{v_0}^2_{L^2(\Omega)} + \eps \norm{u_{e,0}}^2_{L^2(\Omega)}\bigl),
\end{equation}
for  a constant $C = C(\Omega, \omega, T)>0$.

\item If $h$ is   $C^1(\mathbb{R})$ and $h(0)=0$, there exists $\gamma >0$, does not depending on $\varepsilon$, such that,  for every  $(v_0, u_{e, 0}) \in \left(W^{2/3, 6}(\Omega)\cap H^1_0(\Omega)\right)^2$ with $||(v_0,u_{e, 0})||_{W^{2/3, 6}(\Omega)}\leq \gamma $,  there exists a control $f^{\eps} \in L^{6}(\omega \times (0,T))$ such that the solution $(v^{\eps},u^{\eps}_e)$ of  (\ref{Srelaxed-orig}) satisfies:
$$
v^{\eps}(T) = u^\eps_e(T) = 0.
$$
Moreover, the control $f^\varepsilon$ satisfies the estimate:
\begin{equation}\label{bound-control-1}
\norm{f^\varepsilon1_{\omega}}_{L^{6}(Q)}^2 \leq C\bigl( \norm{v_0}_{L^2(\Omega)}^2 + \eps \norm{u_{e,0}}_{L^2(\Omega)}^2\bigl),
\end{equation}
for a constant   $C = C(\Omega, \omega,  T)>0$.
\end{enumerate}
\end{theorem}
The proof of Theorem \ref{controllability4}  is achieved through fixed point arguments, and it will be done in section  \ref{section-nonlinear}. 

\begin{remark} 
In this paper we restrict the dimension to $N=2, 3$, because the bidomain model makes sense only in such dimensions.  Nevertheless, from the mathematical point of view, systems \eqref{Sequiv1}, \eqref{Srelaxed-orig} and \eqref{Srelaxed-orig-1}  make sense for any $N \in \mathbb{N}$ (the 1-d case corresponding to the cable equation) and, taking the initial data in the appropriate space, all the results of this paper can be extended to higher dimensions.
\end{remark}



\section{Carleman inequality}\label{CI}

In this section we  prove Theorem \ref{thm-carleman-syst}.

To simplify the notation, we neglect the index $\eps$ and,  since the only constant which matters in the analysis is $\varepsilon$, we assume that all the other constants are normalized  to be the unity. In this case, the adjoint system (\ref{generaladjoint-1})  reads:

\begin{equation}
	\label{S1-eps-reg-adjoint-concise}
	\begin{cases}
		 -\pt \varphi -\Div (M_e(x)\Grad \varphi) 
		+a(x,t)\varphi= \Div (M_i(x)\Grad \varphi_e) &\qquad \text{ in } Q,\\ 
		 - \eps \pt \varphi_{e}- \Div (M(x)\Grad \varphi_e)=0&\qquad \text{ in } Q, \\
		 \varphi = \varphi_e = 0 &\qquad \text{ on } \Sigma, \\
		\varphi(T) = \varphi_T, \ \varphi_e(T) = \varphi_{e,T}&\qquad \text{ in } \Omega.
	\end{cases}	
\end{equation}

We notice that, if  $\varphi_T$ and $\varphi_{e,T}$ are regular enough, taking $\rho(x,t) = \Div (M_i(x)\nabla \varphi_e(x,t))$,  the pair $(\varphi,\rho)$ satisfies:
\begin{equation}
	\label{S1-eps-reg-adjoint-concise-1}
	\begin{cases}
		 -\pt \varphi -\Div (M_e(x)\Grad \varphi) 
		+a(x,t)\varphi= \rho &\qquad \text{ in } Q,\\ 
		 - \eps \pt \rho- \Div (M(x)\Grad \rho)=0&\qquad \text{ in } Q, \\
		\varphi= \rho =0 &\qquad \text{ on } \Sigma, \\
		   \varphi(T)=\varphi_T, \quad  \rho(T)=\rho_T &\qquad \text{ in } \Omega.
	\end{cases}	
\end{equation}
We prove the Carleman inequality (\ref{observability1112}) using system  \eqref{S1-eps-reg-adjoint-concise-1}.

Before starting the proof of the Carleman inequality, let us first define several weight functions which will be usefull in the sequel.
\begin{lemma} \label{funcaopsi} Let $ \omega_0 $ be  an arbitrary nonempty open set such that $\overline{\omega_0} \subset \omega   \subset \Omega$. There exists a function  
$ \psi \in C^2(\overline{\Omega})$ such that
$$ \psi(x) > 0, \forall x  \in \Omega , \  \  \psi \equiv 0 \  \mbox{on} \ \partial \Omega, \   \   \abs{\nabla \psi(x)} > 0 \ \forall x \in \overline{\Omega \backslash \omega_0}.$$
\end{lemma}
\begin{proof} See \cite{16}.
\end{proof}
Using  Lemma $\ref{funcaopsi}$,  we introduce the weight  functions
\begin{equation}\label{weightfunctions1}
\phi(x,t) = \frac{e^{\lambda(\psi(x) + m ||\psi||)}}{t(T-t)};  \ \  \phi^*(t) = \min_{x \in \overline{\Omega}} \phi(x,t) = \frac{e^{\lambda  m ||\psi||}}{t(T-t)};  
\end{equation}
\begin{equation}\label{weightfunctions2}
\alpha(x,t) = \frac{e^{\lambda(\psi(x) + m||\psi||)} - e^{2\lambda m \norm{\psi}}}{t(T-t)};  \ \ \alpha^*(t) = \max_{x \in \overline{\Omega}}\alpha (x,t)
= \frac{e^{\lambda (m+1) \norm{\psi}} - e^{2\lambda m \norm{\psi}}}{t(T-t)},
\end{equation}
 for a parameter $ \lambda > 0 $ and a constant $m >1$. 
 Here,
$$
 \norm{\psi(x)} = \max_{x \in  \overline{\Omega}}  \abs{\psi(x)}.
$$
\begin{remark}\label{obs2}
From the definition of $\alpha$ and $\alpha^*$ it follows that, for $\lambda$ large enough,  $3\alpha^* \leq 2\alpha$. Moreover,
\begin{equation*}
\phi^*(t) \leq \phi (x,t) \leq e^{\lambda  \norm{\psi}}\phi^*(x,t)
\end{equation*}
and
\begin{equation*}
\abs{\pt \alpha^*} \leq e^{2\lambda  \norm{\psi}}T\phi^2.
\end{equation*}
\end{remark}
\begin{proof}[Proof of Theorem \ref{thm-carleman-syst} ]
For a better comprehension, we divide the proof into several steps.\\

\noindent   \textit{ \textit{Step 1.}} \textit{ First  estimate for the parabolic system.}

In this step we obtain a first Carleman estimate for the adjoint system (\ref{generaladjoint-1}). For that,  we will apply sharp Carleman inequalities, with respect to $\varepsilon$, to the system and get a global estimate of $\varphi$ and $\rho$ in terms of a local integral of $\varphi$ and another in $\rho$.

We consider a set $\omega_1$ such that $\omega_0 \subset \subset \omega_1 \subset \subset \omega$ and  apply the sharp Carleman inequality (\ref{epsilonCarleman}), with $\eps = 1$, and  (\ref{ineq8}) to $\varphi$  and $\rho$, respectively.  We get
\begin{equation} \label{Cestimate1}\begin{split}
&\int \! \! \! \int_{Q} s^{-1}\phi^{-1}e^{2s\alpha}\abs{\varphi_t}^2 dxdt 
+ s^{-1}\int \! \! \! \int_{Q}\phi^{-1}e^{2s\alpha} \sum_{i,j= 1}^N \left | \pxij \varphi \right | ^2dxdt \\ 	
&+ s^3\lambda^4\int \! \! \! \int_{Q}\phi^3e^{2s\alpha}\abs{\varphi}^2dxdt   
+   s\lambda^2\int \! \! \! \int_{Q}\phi e^{2s\alpha} \abs{\nabla \varphi}^2dxdt  \\
 &\qquad  \leq C\biggl( \int \! \! \! \int_{Q} e^{2s\alpha} (\abs{\rho}^2 + \abs{\varphi}^2) dxdt
 + s^3\lambda^4\int \! \! \! \int_{Q_{\omega_1}} \phi^3e^{2s\alpha}\abs{\varphi}^2dxdt    \biggl)
 \end{split} \end{equation}
and 
\begin{equation} \label{ineq777}\begin{split}
&\int \! \! \! \int_{Q} e^{2s\alpha}\abs{\pt \rho}^2 dxdt   + \eps^{-2}\int \! \! \! \int_{Q}e^{2s\alpha}
\sum_{i,j= 1}^N \left | \pxij \rho  \right |^2dxdt  \\  &
\qquad + s^4\lambda^4\eps^{-2}
\int \! \! \! \int_{Q}\phi^4e^{2s\alpha}\abs{\rho}^2dxdt   
+   s^2\lambda^2\eps^{-2}\int \! \! \! \int_{Q}\phi^2e^{2s\alpha} \abs{\nabla \rho}^2dxdt  \\
 & \qquad \qquad \qquad  \leq C e^{\lambda\norm{\psi}} s^4\lambda^4\eps^{-2}\int \! \! \! \int_{Q_{\omega_1}} \phi^4e^{2s\alpha}\abs{\rho}^2dxdt,
\end{split}
\end{equation}
for  $s \geq (T +(1+\norm{a}_{L^{\infty}}^{2/3}) T^2)s_0$ and $\lambda \geq \lambda_0$.
 
Adding (\ref{Cestimate1}) and  (\ref{ineq777}), and absorbing the lower order terms in the right-hand side, we get 
\begin{equation}\label{Cestimate3} \begin{split}
&\int \! \! \! \int_{Q} \phi^{-1}e^{2s\alpha}\abs{\varphi_t}^2 dxdt  
+ \int \! \! \! \int_{Q}\phi^{-1}e^{2s\alpha}  \sum_{i,j= 1}^N  \left| \pxij \varphi \right |^2dxdt  \\ 
&+ s^4\lambda^4\int \! \! \! \int_{Q}\phi^3e^{2s\alpha}\abs{\varphi}^2dxdt   
+   s^2\lambda^2\int \! \! \! \int_{Q}\phi e^{2s\alpha} \abs{\nabla \varphi}^2dxdt  \\
&+\eps^{2}\int \! \! \! \int_{Q} e^{2s\alpha}\abs{\pt \rho}^2 dxdt   
+ \int \! \! \! \int_{Q} e^{2s\alpha} \sum_{i,j= 1}^N \left | \pxij \rho  \right |^2dxdt \\  
&+ s^4\lambda^4\int \! \! \! \int_{Q}\phi^4e^{2s\alpha}\abs{\rho}^2dxdt   
+   s^2\lambda^2\int \! \! \! \int_{Q}\phi^2e^{2s\alpha} \abs{\nabla \rho}^2dxdt   \\
 & \qquad \leq C\biggl(  e^{\lambda||\psi||} s^4\lambda^4\int \! \! \! \int_{Q_{\omega_1}} \phi^4e^{2s\alpha}\abs{\rho}^2dxdt
 + s^4\lambda^4\int \! \! \! \int_{Q_{\omega_1}} \phi^3e^{2s\alpha}\abs{\varphi}^2dxdt    \biggl), 
 \end{split}\end{equation}
 for  $s \geq (T +(1+\norm{a}_{L^{\infty}}^{2/3}) T^2)s_0$ and $\lambda \geq \lambda_0$.
 
\begin{remark}
If we were trying to drive the solution of (\ref{Srelaxed-orig-1}) to zero by means of controls on both equations, inequality (\ref{Cestimate3}) would be sufficient.
\end{remark}
\noindent  \textit{Step 2.} \textit{Estimate of the local integral of $\rho$.}

In this step we estimate the local integral on $\rho$ in the right-hand side of (\ref{Cestimate3}).
This will be done using equation $(\ref{S1-eps-reg-adjoint-concise-1})_1$. Indeed, we consider a function $\xi$ satisfying 
$$
\xi \in C^{\infty}_0(\omega), \ 0 \leq \xi \leq 1, \ \xi(x) = 1 \ \forall x \in \omega_1
$$
and  write
 \begin{equation}\label{ estimate 11 }\begin{split}
Ce^{\lambda||\psi||}s^4\lambda^4\int \! \! \! \int_{Q_{\omega}} e^{2s\alpha}\phi^{4}|\rho|^2\xi dxdt  &= Ce^{\lambda||\psi||}s^4\lambda^4\int \! \! \! \int_{Q_{\omega}} e^{2s\alpha}\phi^{4}\rho
 ( -\varphi_t -\Div (M_e\nabla \varphi )+a\varphi   ) \xi dxdt\\
 &\qquad \qquad :=  E+ F +G.
 \end{split}\end{equation}
In the sequel, we estimate each parcel in the expression above. 
First, we have
\begin{equation*}\begin{split}
E &=Ce^{\lambda||\psi||}s^4\lambda^4\int \! \! \! \int_{Q_{\omega}} s\pt \alpha e^{2s\alpha}\phi^{4}\rho\varphi \xi dxdt \\
& \qquad +Ce^{\lambda||\psi||}s^4\lambda^4\int \! \! \! \int_{Q_{\omega}} e^{2s\alpha}\phi^{3}\phi_t\rho\varphi \xi dxdt 
+ Ce^{\lambda||\psi||}s^4\lambda^4\int \! \! \! \int_{Q_{\omega}} e^{2s\alpha}\phi^{4}\pt \rho \varphi \xi dxdt \\
& :=  E_1 + E_2+ E_3,
\end{split} \end{equation*}
and it is not difficult  to see that
\begin{align}
E_1 + E_2
\leq \frac{1}{10}s^4\lambda^4\int \! \! \! \int_{Q_{\omega}} e^{2s\alpha}\phi^{4}|\rho|^2 dxdt 
+ Ce^{2\lambda||\psi||}s^8\lambda^4\int \! \! \! \int_{Q_{\omega}} e^{2s\alpha}\phi^{8}|\varphi|^2dxdt 
\end{align}
and 
\begin{equation*}\begin{split}
 E_3 \leq \frac{\varepsilon^2}{2}\int \! \! \! \int_{Q_{\omega}} e^{2s\alpha}|\pt \rho|^2 dxdt 
  + Ce^{2\lambda||\psi||}\varepsilon^{-2}s^8\lambda^8\int \! \! \! \int_{Q_{\omega}} e^{2s\alpha}\phi^{8}|\varphi|^2dxdt.
\end{split}\end{equation*}
Next, integrating by parts, we get
\begin{equation*}\begin{split}
e^{-\lambda||\psi||}s^{-4}\lambda^{-4}F &= \sum_{i,j=1}^N\int \! \! \! \int_{Q_{\omega}} s\partial_{x_i} \alpha\,e^{2s\alpha}
\phi^{4}\rho (M_e^{ij} \partial_{x_j}\varphi )\xi dxdt 
 \\&\qquad +\sum_{i,j=1}^N \int \! \! \! \int_{Q_{\omega}} e^{2s\alpha}\phi^{3}
\partial_{x_i}\phi\, \rho  (M_e^{ij} \partial_{x_j}\varphi ) \xi dxdt
\\&\qquad + \sum_{i,j=1}^N\int \! \! \! \int_{Q_{\omega}}
e^{2s\alpha}\phi^{4}\partial_{x_i}\rho\,  (M_e^{ij} \partial_{x_j}\varphi ) \xi  dxdt \\
&\qquad + \sum_{i,j=1}^N \int \! \! \! \int_{Q_{\omega}} e^{2s\alpha}\phi^{4} \rho  (M_e^{ij} \partial_{x_j}\varphi ) \partial_{x_i} \xi  dxdt \\
\end{split}\end{equation*}
and we can show  that 
\begin{equation*}\begin{split}
F  &\leq \frac{1}{10}s^4\lambda^4\int \! \! \! \int_{Q_{\omega}}  e^{2s\alpha}\phi^{4}|\rho|^2dxdt  + \frac{1}{6}s^2\lambda^2\int \! \! \! \int_{Q_{\omega}}  e^{2s\alpha}\phi^2| \nabla \rho|^2dxdt \\
&\qquad  +Ce^{2\lambda||\psi||}s^8\lambda^8\int \! \! \! \int_{Q_{\omega}}  e^{2s\alpha}\phi^{8}|\varphi|^2dxdt + \frac{1}{2}  \int \! \! \! \int_{Q_{\omega}}  e^{2s\alpha}\sum_{i,j=1}^N \left |\pxij \rho \right | ^2dxdt.
 \end{split}\end{equation*} 
Finally, we have 
\begin{equation*}\begin{split}
G \leq &\frac{1}{10}s^4\lambda^4\int \! \! \! \int_{Q_{\omega}} e^{2s\alpha}\phi^{4}|\rho|^2  dxdt \\
& \qquad \qquad +Ce^{2\lambda||\psi||}s^4\lambda^4||a||^2_{L^{\infty}}\int \! \! \! \int_{Q_{\omega}} e^{2s\alpha}\phi^{4} |\varphi|^2    dxdt.
 \end{split}\end{equation*}
Putting  $E$, $F$ and $G$ together in (\ref{Cestimate3}), we obtain
\begin{equation}\label{Cestimate4}\begin{split}
&\int \! \! \! \int_{Q} e^{2s\alpha}|\varphi_t|^2 dxdt  + \int \! \! \! \int_{Q}e^{2s\alpha} \sum_{i,j= 1}^N
\left | \pxij \varphi \right |^2dxdt 
\\&\qquad + s^4\lambda^4\int \! \! \! \int_{Q}\phi^4e^{2s\alpha}|\varphi|^2dxdt  +   s^2\lambda^2\int \! \! \! \int_{Q}\phi^2 e^{2s\alpha} |\nabla \varphi|^2dxdt 
\eps^2\int \! \! \! \int_{Q} e^{2s\alpha}|\pt \rho|^2 dxdt  
\\&\qquad \qquad+ \int \! \! \! \int_{Q}e^{2s\alpha} \sum_{i,j= 1}^N \left | \pxij \rho \right |^2dxdt 
+ s^4\lambda^4\int \! \! \! \int_{Q}\phi^4e^{2s\alpha}|\rho|^2dxdt   
\\&\qquad \qquad \qquad +   s^2\lambda^2\int \! \! \! \int_{Q}\phi^2 e^{2s\alpha} |\nabla \rho|^2dxdt   \\  
&\leq C e^{2\lambda||\psi||}\varepsilon^{-2}s^8\lambda^8\int \! \! \! \int_{Q_{\omega}} e^{2s\alpha}\phi^{8}|\varphi|^2 dxdt,
 \end{split}\end{equation}
for  $s \geq (T +(1+\norm{a}_{L^{\infty}}^{2/3}) T^2)s_0$ and $\lambda \geq \lambda_0$.

Using (\ref{Cestimate4}), we can prove that, for every $\eps >0$, system (\ref{Srelaxed-orig-1})  is null controllable. However,  the sequence of controls obtained in this way will not be bounded when $\eps \rightarrow 0^+$.  Therefore, we need to go a step further and improve estimate (\ref{Cestimate4}). This is the goal of the next step.

\vskip0.3cm

\noindent \textit{Step 3.} \textit{Weighted energy inequality.}\\
The reason why we do not get a bounded sequence of controls out of step 2 is because of the term $\eps^{-2}$ in the right-hand side of (\ref{Cestimate4}).   In this step  we  prove  a weighted energy inequality for equation $(\ref{S1-eps-reg-adjoint-concise-1})_2$, which will be used to  compensate this $\eps^{-2}$ term.

Let us introduce the function 
$$y = e^{\frac{3}{2}s\alpha^*}\rho.$$
 This new function satisfies
\begin{equation}\label{eq222} 
\begin{cases}
\eps \partial_ty -\Div ( M(x)\nabla y) = \eps\frac{3}{2}s\pt \alpha^* e^{\frac{3}{2}s\alpha^*}\rho  &\qquad \mbox{ in } Q, \\
y = 0 &\qquad \mbox{ on } \Sigma, \\
y(0) = y(T) = 0 &\qquad \mbox{ in }  \Omega.
\end{cases}
\end{equation}
Multiplying   (\ref{eq222}) by $y$ and integrating over $\Omega$, we get 
$$
\frac{\varepsilon}{2}\frac{d}{dt}||y(t)||^2_{L^2(\Omega)} +C||\nabla y(t)||^2_{L^2(\Omega)}  \leq \varepsilon\frac{3}{2}\int_{\Omega}s\pt \alpha^*(t) e^{\frac{3}{2}s\alpha^*(t)}\rho(t) y(t) dx.
$$
Integrating this last inequality form $0$ to $T$ and using Poincar\'{e}'s and Young's  inequalities,  it is not difficult to see  that
\begin{equation}\label{wenergy}
\iint_{Q}  e^{3s\alpha^*}|\rho|^2dxdt   \leq C\varepsilon^2 e^{4\lambda||\psi||} \iint_{Q}s^4\phi^4 e^{2s\alpha}|\rho|^2  dxdt. 
\end{equation}
Finally, from  (\ref{Cestimate4}) and  (\ref{wenergy}),  we obtain
\begin{equation}\label{equat123}
\iint_{Q}  e^{3s\alpha^*}|\rho|^2dxdt   \leq C e^{6\lambda||\psi||}s^8\lambda^4 \iint_{Q_{\omega}}\phi^8 e^{2s\alpha}|\varphi|^2  dxdt.
\end{equation}

This last  estimate gives a global estimate of $\rho$ in terms of a local integral of $\phi$, with a constant $C$ which is bounded with respect to $\varepsilon$.

\vskip0.2cm

\noindent  \textit{Step 4.} \textit{Last estimates  and conclusion.}\\
In order to finish the proof of Theorem \ref{thm-carleman-syst},  we combine inequality (\ref{equat123}) and a slightly different Carleman inequality  to the equation $(\ref{S1-eps-reg-adjoint-concise-1})_1$. Indeed, the following Carleman inequality holds: 
\begin{equation}\label{Cestimate1323}\begin{split}
&\int \! \! \! \int_{Q} s^{-1}\phi^{-1}e^{3s\alpha}|\varphi_t|^2 dxdt  + s^{-1}\int \! \! \! \int_{Q}\phi^{-1}e^{3s\alpha} \sum_{i,j= 1}^N \left |\pxij \varphi \right |^2dxdt  \\ 
&\qquad \qquad + s^3\lambda^4\int \! \! \! \int_{Q}\phi^3e^{3s\alpha}|\varphi|^2dxdt   
+   s\lambda^2\int \! \! \! \int_{Q}\phi e^{3s\alpha} |\nabla \varphi|^2dxdt   \\
 &\qquad  \leq C\biggl( \int \! \! \! \int_{Q} e^{3s\alpha}|\rho|^2 dxdt+ s^3\lambda^4\int \! \! \! \int_{Q_{\omega}} \phi^3e^{3s\alpha}|\varphi|^2dxdt    \biggl),
\end{split}\end{equation}
for  $s \geq (T +(1+\norm{a}_{L^{\infty}}^{2/3}) T^2)s_0$ and $\lambda \geq \lambda_0$, where $\varphi$ is, together with $\rho$, solution of  (\ref{S1-eps-reg-adjoint-concise-1}).

Notice that here we have  just changed the weight  $e^{2s\alpha}$ by $e^{3s\alpha}$. The proof of (\ref{Cestimate1323}) is exactly the same as  the proof  of Theorem \ref{carleman-calor},  just taking the appropriate change of variable in (\ref{change-carleman}).\

Next, since $e^{3s\alpha} \leq e^{3s\alpha^*}$, we have
$$
 \int \! \! \! \int_{Q}e^{3s\alpha}|\rho|^2dxdt \leq  \int \! \! \! \int_{Q}e^{3s\alpha^*}|\rho|^2dxdt
$$
and by (\ref{equat123}), we have that 
$$
\int \! \! \! \int_{Q}e^{3s\alpha}|\rho|^2dxdt \leq C e^{6\lambda||\psi||}s^8\lambda^4 \iint_{Q_{\omega}}\phi^8 e^{2s\alpha}|\varphi|^2  dxdt .
$$
From (\ref{equat123})  and (\ref{Cestimate1323}), it follows that
 \begin{equation}\label{observability1}
 \int \! \! \! \int_{Q}e^{3s\alpha}|\rho|^2dxdt +s^3\lambda^4\int \! \! \! \int_{Q}\phi^3e^{3s\alpha}|\varphi|^2dxdt   
 \leq C e^{6\lambda||\psi||}s^8\lambda^4 \iint_{Q_{\omega}}\phi^8 e^{2s\alpha}|\varphi|^2  dxdt,
\end{equation}
which is exactly (\ref{observability1112}). 

By density, we can show that  (\ref{observability1}) remains true when we  consider  initial data in $L^2(\Omega)$. Therefore,  the proof of Theorem \ref{thm-carleman-syst} is finished.

\end{proof}


\section{Null controllability for the linearized system}\label{control-linearsystem}
This section is devoted to proving the null controllability of linearized equation (\ref{Srelaxed-orig-1}). It will be done by showing the observability inequality (\ref{observability-inequality}) for  the adjoint system (\ref{generaladjoint-1}), and solving a minimization problem.
The arguments used here  are  classical in control theory for linear PDE's. Hence, we just give a sketch of the proof. 

\begin{proof}[Proof of Theorem \ref{controllability2}]

Combining the standard energy inequalities for system (\ref{S1-eps-reg-adjoint-concise-1}) and the Carleman inequality given by Theorem \ref{thm-carleman-syst}, we can show the following observability inequality for the solutions of (\ref{S1-eps-reg-adjoint-concise-1}):
\begin{equation}\label{observability2}
||\varphi(0)||^2_{L^2(\Omega)}+\varepsilon ||\rho(0)||_{L^2(\Omega)}^2  \leq   e^{C(1 + 1/T + ||a||^{2/3}_{L^\infty} + ||a||_{L^\infty}T)} \iint_{Q_{\omega}}|\varphi|^2  dxdt,
\end{equation}
where $C = C(\Omega, \omega) $ is a positive constant.

Next, since $\rho(x,t) = \Div (M(x)\nabla \varphi_e(x,t))$ and $\varphi_e = 0$ on $\partial \Omega$, we have that
$$
||\varphi_e(t)||_{H^2(\Omega)} \leq C||\rho(t)||_{L^2(\Omega)},
$$
for all $t \in [0,T]$. Therefore,  it follows from (\ref{observability2})  that 
\begin{equation}\label{observability3}
||\varphi(0)||^2_{L^2(\Omega)} +\varepsilon ||\varphi_e(0)||_{L^2(\Omega)}^2 \leq  e^{C(1 + 1/T + ||a||^{2/3}_{L^\infty} + ||a||_{L^\infty}T)} \iint_{Q_{\omega}}|\varphi|^2  dxdt,
\end{equation}
which is the observability inequality (\ref{observability-inequality}).

From (\ref{observability3}) and the density of smooth solutions in the space  of solutions of (\ref{S1-eps-reg-adjoint-concise})  with initial data in $L^2(\Omega)$, we see that the above observability inequality is satisfied by all solutions of  (\ref{generaladjoint-1})  with initial data in $L^2(\Omega)$. 

Now, in order to obtain the null controllability for  linear system (\ref{Srelaxed-orig-1}), we solve, for any $\delta >0$,  the following  minimization problem: 

\begin{equation} \label{Jfunctional}
\begin{split}
&\text{Minimize}\, \J_\delta(\varphi_T,\varphi_{e,T}), \text{ with }\\
&\J_\delta(\varphi_T,\varphi_{e,T})=\Biggl\{\frac{1}{2}\int_0^T\int_{\omega} \abs{\varphi^{\eps}}^2\dx\dt
+\eps(u_{e,0}, \varphi^{\eps}_e(0))\\
&\qquad \qquad \qquad \qquad \qquad \qquad
+(v_0,\varphi^{\eps}(0))+\delta (\norm{\varphi_T}_{L^2(\Omega)}+\varepsilon^{1/2} \norm{\varphi_{e,T}}_{L^2(\Omega)})\Biggl\},
\end{split}
\end{equation}
where $(\varphi, \varphi_e)$ is the solution of the adjoint problem (\ref{generaladjoint-1}) with initital data $(\varphi_T,\varphi_{e,T}) \in L^2(\Omega)^2$. 

It is an easy matter to check  that $\J_\delta$ is strictly convex and  continuous. So, in order to guarantee the existence of a minimizer, 
the only thing remaining  to prove is the coercivity of $\J_{\delta}$.

Using the observability inequality (\ref{observability-inequality}) for the adjoint system  (\ref{generaladjoint-1}),
the coercivity of $\J_{\delta}$ is straightfoward.
Therefore, for each $\delta >0$, there exists a unique minimizer  $(\varphi_{e,T}^{\delta},\varphi_T^{\delta})$ of $\J_\delta$. 
Let us denote by $ \varphi^{\eps,\delta}$ the corresponding solution to  (\ref{generaladjoint-1}) associated to this minimizer. Taking $f^{\eps, \delta} = \varphi^{\eps,\delta}1_{\omega}$ as a control for (\ref{Srelaxed-orig-1}), the duality between (\ref{Srelaxed-orig-1}) and (\ref{generaladjoint-1}) gives
\begin{equation}\label{aproximatecontrollability}
||v^{\varepsilon, \delta}(T)||_{L^2(\Omega)} + \varepsilon^{1/2}||u^{\varepsilon, \delta}_e(T)||_{L^2(\Omega)} \leq \delta,
\end{equation}
where $(v^{\eps, \delta},u_e^{\eps, \delta}) $ is the solution of (\ref{Srelaxed-orig-1})  associated to the control $f^{\eps, \delta}$.
It also gives
\begin{equation}\label{bound1}
||f^{\eps, \delta}1_{\omega}||^2_{L^2(Q)} \leq C\bigl(  ||v_0||^2_{L^2(\Omega)} + \eps||u_{e,0}||^2_{L^2(\Omega)} \bigl).
\end{equation}
From (\ref{aproximatecontrollability}) and (\ref{bound1}), we get a control $f^\eps$
(the weak limit of a subsequence of $f^{\eps, \delta}1_{\omega}$ in $L^2(\omega\times (0,T))$)
that drives the solution of  (\ref{Srelaxed-orig-1})  to zero at time $T$. 
From (\ref{bound1}), we have the following estimate on the control $f^\eps$, 
\begin{equation}\label{bound3}
||f^\eps1_{\omega}||^2_{L^2(Q)} \leq C\bigl( ||v^\eps_0||^2_{L^2(\Omega)} + \eps||u^\eps_{e,0}||^2_{L^2(\Omega)} \bigl).
\end{equation}
This finishes the proof of Theorem \ref{controllability2}.
\end{proof}


\section{The nonlinear system} \label{section-nonlinear}

In this section we prove Theorem \ref{controllability4}.  The proof is achieved through fixed point arguments. 

 
\textit{Proof of Theorem \ref{controllability4} (case 1):} We consider the following linearization of  system (\ref{Srelaxed-orig}):
\begin{equation}
	\label{Srelaxed-orig-12}
	\begin{cases}	
	 c_{m} \pt v^{\eps}-\frac{\mu}{\mu+1}\Div\bigl(\bM_{e}(x)\Grad v^{\eps} \bigr)+
g(z)v^{\eps} = f^{\eps}1_{\omega}&\qquad \text{ in } Q, \\
 \eps \pt u^{\eps}_e-\Div\bigl(\bM(x)\Grad u^{\eps}_{e}\bigr)  = \Div\bigl(\bM_{i}(x)\Grad v^{\eps}\bigr)&\qquad \text{ in } Q, \\
 v^{\eps} =  u^{\eps}_e = 0  &\qquad \text{ on } \Sigma, \\
 v^{\eps}(0) = v_0, \ u^{\eps}_e(0) = u_{e,0} &\qquad \text{ in } \Omega,
	\end{cases}	
\end{equation}
where 
\begin{equation}\label{linearization-h}
g(s)=  
\begin{cases}
\displaystyle \frac{h(s)}{s}, \ \text{if} \ |s| >0, \\
h'(0), \ \text{if} \ s= 0.
\end{cases}
\end{equation}
It follows from Theorem \ref{controllability2} that, for each $(v_0,  u_{e,0})  \in L^2({\Omega})^{2}$ and $z \in L^{2}(Q)$,  there exists a control function $f^{\eps} \in L^2(Q)$ such that the solution of (\ref{Srelaxed-orig-12}) satisfies:
$$
v^{\eps}(T) = u^{\eps}_e(T) = 0.
$$
As we said before, the idea is to use a fixed point argument. For that, we need the following  generalized version  of Kakutani's fixed point Theorem,  due to Glicksberg \cite{Glick}. 
\begin{theorem}\label{Glicksberg}
Let $B$ be a non-empty convex, compact subset of a locally convex topological vector space $X$. If $\Lambda : B \longrightarrow B$  is a convex set-valued  mapping with closed graph and $\Lambda(B)$ is closed, then $\Lambda$ has a fixed point. 
\end{theorem}
In order to apply Glicksberg`s Theorem, we define a mapping  $\Lambda : B \longrightarrow X$ as follows
\begin{equation*} \begin{split}
&\Lambda (z) = \{ v^{\eps}; \  (v^{\eps}, u_e^{\eps}) \ \text{is a solution of (\ref{Srelaxed-orig-12}), such that} \ v^{\eps}(T) = u_e^{\eps}(T) = 0, \\	
&\qquad \qquad\qquad\qquad \qquad\qquad \qquad\qquad
\text{for a control } \  f^{\eps}  \ \text{satisfying} \ (\ref{bound-control})\}.
\end{split}\end{equation*}
Here, $X = L^2(Q)$ and $B$ is the ball
\begin{equation*} \begin{split}
&B = \{ z \in L^2(0,T, H^1_0(\Omega)), \pt z \in  L^2(0,T, H^{-1}(\Omega)); \\ 
&\qquad \qquad  \qquad \qquad \qquad \qquad  ||z||^2_{L^2(0,T;H^1_0(\Omega))}  + ||\pt z||^2_{L^2(0,T;H^{-1}(\Omega))} \leq M  \}.
\end{split}\end{equation*}
It is easy to see that $\Lambda$ is well defined and that $B$ is a  convex and compact subset of $L^2(Q)$.

Let us now prove  that $\Lambda$ is convex, compact and has closed graph.
\vskip0.1cm

$\bullet$ $\Lambda(B) \subset B$. 

Let $z \in B$ and $v^{\eps} \in \Lambda(z)$. Since $v^{\eps}$ satisfies $(\ref{Srelaxed-orig-12})_1$, the following inequality holds
\begin{equation}
||v^{\eps}||^2_{L^2(0,T;H^1_0(\Omega))} + ||\pt v^{\eps}||^2_{L^2(0,T;H^{-1}(\Omega))} \leq K_1.
\end{equation}
In this way,  if $z \in B$ then $\Lambda (z) \subset B$, if we take $M = K_1$.\\

$\bullet$  $\Lambda(z)$ is closed in $L^2(Q)$.

Let $z \in B$ fixed, and $ v_n^{\eps} \in \Lambda(z)$, such that $ v_n^{\eps} \rightarrow v^{\eps}$. Let us prove that $v^{\eps} \in \Lambda(z)$.

\noindent In fact, by definition we have that $v_n^{\eps}$ is, together with a function $u_{e,n}^{\eps}$, and a control $f_n^{\eps}$, the solution of (\ref{Srelaxed-orig-12}), with $||f_n^{\eps}1_{\omega}||^{2}_{L^2(Q)} \leq C\bigl(  ||v_0||^{2}_{L^2(\Omega)} + \eps ||u_{e,0}||^{2}_{L^2(\Omega)}\bigl)  $. 
Therefore, we can extract a subsequence of $f_n^{\eps}$, denoted by the same index, such that 
$$
f_n^{\eps}1_{\omega} \rightarrow f^{\eps}1_{\omega} \ \text{weakly in } \ L^2(Q).   
$$
Since $f_n^{\eps}$ is bounded, we can argue as in the previous section and show that
\begin{equation}
||v_n^{\eps}||^2_{L^2(0,T;H^1_0(\Omega))} + ||\pt v^{\eps}_{n}||^2_{L^2(0,T;H^{-1}(\Omega))} \leq M.
\end{equation}
Hence,
\begin{equation}
\left |   
\begin{array}{ll}
v_n^{\eps} \rightarrow v^{\eps} \ \text{weakly in } \ L^2(0,T; H^1_0(\Omega)), \nonumber\\
v_n^{\eps} \rightarrow v^{\eps} \ \text{strongly in } \ L^2(Q), \nonumber \\
\pt v_{n}^{\eps} \rightarrow \pt v^{\eps} \ \text{weakly in } \ L^2(0,T; H^{-1}(\Omega)).
\end{array}
\right. 
\end{equation}
Using  the  converges above and $(\ref{Srelaxed-orig-12})_2$,  we see that there exists a function $u_e^{\eps}$ such that 
\begin{equation}
\left |   
\begin{array}{ll}
u_{e,n}^{\eps} \rightarrow u_e^{\eps} \ \text{weakly in } \ L^2(0,T; H^1_0(\Omega)), \nonumber\\
u_{e,n}^{\eps} \rightarrow u_e^{\eps} \ \text{strongly in } \ L^2(Q), \nonumber \\
 \pt u_{e,n}^{\eps} \rightarrow \pt u_{e}^{\eps} \ \text{weakly in } \ L^2(0,T; H^{-1}(\Omega)).
\end{array}
\right. 
\end{equation}
It follows that $(v^{\eps},u_e^{\eps})$ is a controlled solution of  (\ref{Srelaxed-orig-12}) associated to the control $f^{\eps}$.
Hence, $v^{\eps} \in \Lambda(z)$ and $\Lambda(z)$ is closed and compact in $L^2(Q)$.\\

$\bullet$ $\Lambda$ has closed graph in $L^2(Q)\times L^2(Q)$.

We need to prove that if $z_n \rightarrow z $, $v_n^{\eps} \rightarrow v^{\eps}$ strongly in $L^2(Q)$ and
$v_n^{\eps} \in \Lambda(z_n)$, then $v^{\eps} \in \Lambda(z)$. Using the two previous steps, it is easy to show that $v^{\eps} \in \Lambda(z)$.

Therefore, we can apply Glicksberg's Theorem to  conclude that  $\Lambda $ has a fixed point.
This proves Theorem \ref{controllability4} in the case where the nonlinearity is a  $C^1$ globally Lipschitz function.
\\


\textit{Proof of Theorem \ref{controllability4} (case 2):}
We consider the linear system:
\begin{equation}
	\label{linearization-cubic}
	\begin{cases}
		 c_{m} \pt v^{\eps}-\frac{\mu}{\mu+1}\Div\bigl(\bM_{e}(x)\Grad v^{\eps} \bigr)+a(z)v^{\eps} = f^{\eps}1_{\omega}
&\qquad \text{ in } Q, \\
\eps \pt u_e^{\eps}-\Div\bigl(\bM(x)\Grad u^{\eps}_{e}\bigr)  = \Div\bigl(\bM_{i}(x)\Grad v^{\eps}\bigr)
&\qquad \text{ in } Q, \\
v^{\eps} = u^{\eps}_e = 0  &\qquad \text{ on } \Sigma, \\
v^{\eps}(0) = v_0, \ u_e^{\eps}(0) = u_{e,0} &\qquad \text{ in } \Omega,
	\end{cases}	
\end{equation}
with $(v_0, u_{e, 0})\in \left(W^{2/3, 6}(\Omega)\cap H^1_0(\Omega)\right)^2$, $z \in L^{\infty}(Q)$ and 
$$
a(z)= \int_0^1 \frac{d h}{dz}(sz)ds.
$$
Arguing as in the proof Theorem \ref{controllability2}, we can show the null controllability of \eqref{linearization-cubic} with controls in $L^2(\omega\times (0,T))$. However,  these $L^2$ controls are not sufficient to apply fixed point  arguments and obtain the null controllability of the nonlinear system \eqref{Srelaxed-orig}. For this reason, we modify a little the functional (\ref{Jfunctional}), obtaining  controls which  will allow us to employ Schauder's fixed point Theorem. Indeed, for any $\delta >0$, we consider the problem:
\begin{equation} \label{Jfunctional-2}
\begin{split}
&\text{Minimize}\, \J_\delta(\varphi_T,\varphi_{eT}), \text{ with }\\
&\J_\delta(\varphi_T,\varphi_{eT})=\Biggl\{\frac{1}{2}\int_0^T\int_{\omega} e^{2s\alpha}\phi^8\abs{\varphi^{\eps}}^2\dx\dt
+\eps(u_{e,0}, \varphi^{\eps}_e(0))\\
&\qquad \qquad \qquad \qquad \qquad
+(v_0,\varphi^{\eps}(0))+\delta \Bigl(\norm{\varphi_T}_{L^2(\Omega)}+\epsilon^{1/2}\norm{\varphi_{e,T}}_{L^2(\Omega)}\Bigl)\Biggl\},
\end{split}
\end{equation}
where $(\varphi^{\eps}, \varphi_e^{\eps})$ is the solution of the adjoint system (\ref{generaladjoint-1}) with initital data $(\varphi_T,\varphi_{e,T}) \in L^2(\Omega)^2$.

As in section \ref{control-linearsystem}, we show that problem (\ref{Jfunctional-2}) has a unique minimizer $(\varphi^{\eps,\delta}, \varphi_e^{\eps, \delta})$.
Defining $f^{\eps, \delta} = e^{2s\alpha}\phi^8\varphi^{\eps, \delta}$ and using the fact that $\varphi^{\eps, \delta}$ is, together with a $\varphi_e^{\eps, \delta}$,
the solution of (\ref{generaladjoint-1}),  we see that $f^{\eps, \delta}$ is a solution of a parabolic  equation, with homogeneous Dirichlet boundary conditions, null initial data and a right-hand side in $L^2(Q)$. Hence, $f^{\eps, \delta} \in L^2(0,T; H^2(\Omega))\cap H^1(0,T; L^2(\Omega))$ and in particular we have that $f^{\eps, \delta} \in L^6(Q)$. 

Arguing as in \cite[Lemma 5]{Khodja-1}, we can easily show that
\begin{equation}\label{Linftyestimate}
||f^{\eps, \delta}1_{\omega}||_{L^{6}(Q)}^2 \leq C\bigl( \norm{v_0}_{L^2(\Omega)}^2 + \eps \norm{u_{e,0}}_{L^2(\Omega)}^2\bigl), 
\end{equation}
where $C>0$ is independent of $\epsilon$ and $\delta$. 

Moreover, the solution  $(v^{\eps, \delta},u_e^{\eps, \delta})$ of \eqref{linearization-cubic}, associated to  $f^{\eps, \delta}$, satisfies 
\begin{equation}\label{aproximatecontrollability-2}
||v^{\varepsilon, \delta}(T)||_{L^2(\Omega)} + \varepsilon^{1/2}||u^{\varepsilon, \delta}_e(T)||_{L^2(\Omega)} \leq \delta.
\end{equation}

Taking the limit when $\delta \rightarrow 0^+$, we get a control $f^\eps \in L^{6}(Q)$ (the weak limit of a subsequence of $f^{\eps, \delta}$)  such that the associated solution $(v^\eps,u_e^\eps)$ to (\ref{linearization-cubic}) satisfies 
$$
v^\eps(T) = u^\eps_e(T) = 0.
$$

Next, we define a  map $F: L^{\infty}(Q) \longrightarrow L^{\infty}(Q)$ which, to each $z \in L^{\infty}(Q)$, associates $v^{\eps}$ the solution, together with $u_e^{\eps}$, of \eqref{linearization-cubic} corresponding to $z$ and to the control $f^{\eps}$ built above. Note that  this application is well defined  since, from the regularity theory for parabolic equations (see, for instance, \cite{Landyz}), we have that $v^{\eps} \in X := L^6(0,T; W^{2,6}(\Omega)) \cap W^{1,6}(0,T; L^6(\Omega)) \cap C([0,T]; W^{2/3,6}(\Omega)\cap H^1_0(\Omega))$.

Let us now consider the set $A$ defined as 
$$
A := \left\{ z \in L^{\infty}(Q); ||z||_{L^{\infty}(Q)} \leq 1 \right\}. 
$$

It is clear that $A$ is a  convex closed subset of $L^{\infty}(Q)$. From \eqref{Linftyestimate}, which still holds for $f^\eps$,  and the smallness assumption on the initial data, we can easily show  that  $F$ is continuous and that $F(A) \subset A$. Finally, since the space $X$ is compactly embedded in $L^{\infty}(Q)$,  we have that $F(A)$ is compact in $L^{\infty}(Q)$. Therefore, $F$ has a fixed point and the the  proof of case $2$ of Theorem \ref{controllability4} is finished.
 

\section{Appendix: Some technical results}\label{appendix-carleman1}
In this section we prove the two  sharp Carleman inequalities  used in the proof of Theorem \ref{thm-carleman-syst}.

We consider the parabolic equation
\begin{equation}\label{eq1-1}
\begin{cases}
-\pt v(x,t)   - \displaystyle\sum_{i,j =1}^N \pxi(a_{ij}(x)\pxi v(t,x)) = g(x,t) &   \qquad   \mbox{ in }  Q ,  \\
v= 0 & \qquad    \mbox{ on }   \Sigma, \\
v(T)  = v_T       &  \qquad  \mbox{ in }  \Omega,
\end{cases}
\end{equation}
where $v_T \in L^2(\Omega)$ and  $g \in L^2(Q)$. 
 
We assume that  the matrix $a_{ij}$ has the form
$$
a_{ij} = \frac{M_{ij}}{\eps},
$$
and $(M_{ij})_{ij}$ is an elliptic matrix, i.e.,  there exists $\beta >0$ such that $\sum_{i,j}^N M_{ij}\xi_j\xi_i \geq \beta |\xi|^2$ for all $\xi \in \mathbb{R}^N$.

\subsection{A degenerating Carleman inequality}

The first sharp Carleman inequality we prove is the following.
\begin{theorem}\label{carleman-calor}
For any $0< \epsilon \leq 1$, there exist $\lambda_0  = \lambda_0(\Omega, \omega) \geq 1 $ and $s_0 =  s_0(\Omega, \omega) \geq 1 $ such that, for every $ \lambda \geq  \lambda_0$ and   $s \geq s_0 (T + T^2)$,  the solution $v$ of  (\ref{eq1-1}) satisfies 
\begin{equation} \label{epsilonCarleman}\begin{split}
&\int \! \! \! \int_{Q} s^{-1}\phi^{-1}e^{2s\alpha}|\pt v|^2 dxdt  
+ s^{-1}\eps^{-2}\int \! \! \! \int_{Q}\phi^{-1}e^{2s\alpha} \sum_{i,j= 1}^N \left | \pxij v \right | ^2dxdt \\  
&\qquad \qquad + s^3\lambda^4\eps^{-2}\int \! \! \! \int_{Q}\phi^3e^{2s\alpha}|v|^2dxdt   
+   s\lambda^2\eps^{-2}\int \! \! \! \int_{Q}\phi e^{2s\alpha} |\nabla v|^2dxdt     \\
 & \qquad \leq C\biggl(  \int \! \! \! \int_{Q}e^{2s\alpha}|g|^2dxdt+ s^3\lambda^4\eps^{-2}\int \! \! \! \int_{Q_{\omega}} \phi^3e^{2s\alpha}|v|^2dxdt     \biggl), 
 \end{split}\end{equation}
with $C>0$  depending only on $\Omega$, $\omega_0$, $\psi$ and $\beta$.
\end{theorem}
 \begin{proof}
 For $s > 0$ and $\lambda >0$, we consider the change of variable
 \begin{equation}\label{change-carleman}
 w(t,w) = e^{s\alpha} v(t,w),
 \end{equation}
 which implies
 $$
 w(T,x) = w(0,x) = 0.
 $$
 
Using the fact that $v$ is the solution of \eqref{eq1-1}, we write
\begin{equation}\label{w-equation}
L_1w + L_2 w = g_s,
\end{equation}
where
\begin{equation}\label{L1}
L_1w =  - \pt w + 2s\lambda \sum_{i,j = 1}^N \phi a_{ij}\pxj \psi\, \pxi w + 2s\lambda^2 
\sum_{i,j = 1}^N \phi a_{ij}\pxi \psi\,\pxj \psi w,
\end{equation}
\begin{equation}\label{L2}
L_2 w = -\sum_{i,j = 1}^N \pxi (a_{ij}\pxj w ) -s^2 \lambda^2 \sum_{i,j = 1}^N \phi^2 a_{ij}\pxi\psi\,\pxj \psi \,w 
+s\pt \alpha \, w
\end{equation}
and
\begin{equation}
g_s = e^{s\alpha}g + s\lambda^2\sum_{i,j = 1}^N\phi a_{ij}\pxi \psi \, \pxj \psi \, w 
-s\lambda \sum_{i,j =1}^N\phi \pxi (a_{ij}\pxj \psi )\,w.   
\end{equation}
From (\ref{w-equation}), we have that 
\begin{equation} \label{id11}
||L_1w||_{L^2(Q)}^2 + ||L_2w||_{L^2(Q)}^2 + 2 (L_1w,L_2w)_{L^2(Q)} = ||g_s||_{L^2(Q)}^2.
\end{equation}
The rest of the proof is devoted  to analyze the terms appearing in $  (L_1w,L_2w)_{L^2(Q)}$. 

 First, we write 
$$
(L_1w,L_2w)_{L^2(Q)} = \sum_{i,j = 1}^N I_{ij},
$$
where $I_{ij}$  is the inner product in $L^2(Q)$ of the $ith$ term in the expression of $L_1w$ and the $jth$ term in $L_2w$. 

After a long, but straightforward, calculation, we can show that the following estimate holds
 
\begin{align}
2(L_1w,L_2w)_{L^2(Q)}&\geq  2 s^3\lambda^4\beta^2\varepsilon^{-2}\int \! \! \! \int_{Q}\phi^3|\nabla\psi|^4|w|^2dxdt +2s\lambda^2\beta^2\varepsilon^{-2}\int \! \! \! \int_{Q}\phi|\nabla\psi|^2|\nabla w|^2dxdt\nonumber \\
& - C\varepsilon^{-2}\biggl(    T^2s^2 \lambda^4 +    Ts^2\lambda^2 +        T^2s +    s^3\lambda^3   +  T s^2\lambda     \biggl)\int \! \! \! \int_{Q}\phi^3|w|^2dxdt  \nonumber \\  &  -C\varepsilon^{-2}(  s\lambda +\lambda^2)\int \! \! \! \int_{Q}\phi|\nabla w|^2dxdt. \label{112}   
\end{align}

We take $\lambda \geq \lambda_0$ and $s \geq s_0(T +T^2)$,  and it  follows,  from Remark \ref{obs11} below, that 

\begin{align}\label{610}
& 2(L_1w,L_2w)_{L^2(Q)} +  2 s^3\lambda^4\beta^2\varepsilon^{-2}\int \! \! \! \int_{Q_{\omega_0}}\phi^3|w|^2dxdt \nonumber \\
 &+2s\lambda^2\beta^2\varepsilon^{-2}\int \! \! \! \int_{Q_{\omega_0}}\phi|\nabla w|^2dxdt  \nonumber \\
 & \geq  2 s^3\lambda^4\beta^2\varepsilon^{-2}\int \! \! \! \int_{Q}\phi^3|w|^2dxdt +2s\lambda^2\beta^2\varepsilon^{-2}\int \! \! \! \int_{Q}\phi|\nabla w|^2dxdt.
\end{align}

\begin{remark}\label{obs11}
 Since $\overline{\Omega \backslash \omega_0}$ is compact and $|\nabla \psi| > 0  \ \mbox{on} 
 \ \overline{\Omega \backslash \omega_0}$, there exists $\delta >0$ such that 
 $$
 \beta |\nabla\psi| \geq \delta \ \mbox{on}  \ \overline{\Omega \backslash \omega_0}.
 $$ 
  \end{remark}

Putting (\ref{610})  in  (\ref{id11}), we get
\begin{equation}\begin{split}\label{eq31}
& ||L_1w||_{L^2(Q)}^2 + ||L_2w||_{L^2(Q)}^2 + 2\beta^{-2} s^3\lambda^4\delta^4\varepsilon^{-2}\int \! \! \! \int_{Q}\phi^3|w|^2dxdt 
 \\&\qquad \qquad \qquad+   2s\lambda^2\delta^2\varepsilon^{-2}\int \! \! \! \int_{Q}\phi |\nabla w|^2dxdt  \\ 
&\leq   ||g_s||_{L^2(Q)}^2  + 2\beta^{-2} s^3\lambda^4\delta^4\varepsilon^{-2}\int \! \! \! \int_{Q_{\omega_0}} \phi^3|w|^2dxdt \\ 
& \qquad \qquad \qquad+   2 s\lambda^2\delta^2\varepsilon^{-2}\int \! \! \! \int_{Q_{\omega_0}}\phi|\nabla w|^2dxdt.  
\end{split}\end{equation}

Now we deal with the local integral involving $\nabla w$ on the right-hand side of (\ref{eq31}).
To this end, we introduce a cutt-off function $\xi$ such that
$$
\xi \in C^{\infty}_0(\omega), \ 0 \leq \xi \leq 1, \ \xi(x) = 1 \ \forall x \in \omega_0.
$$
Using the ellipticity condition on $a_{ij}$, we can prove that
\begin{equation*}\begin{split} 
&\beta\varepsilon^{-1}\int \! \! \! \int_{Q_{\omega}}\phi \xi^2|\nabla w|^2 dxdt 
\\  &\qquad \leq  C \biggl(\int \! \! \! \int_{Q} L_2w\phi\xi^2wdxdt +  (sT + \varepsilon^{-1}s^2 \lambda^2) \int \! \! \! \int_{Q_{\omega}}   \phi^3|w|^2dxdt  \nonumber  \\
&\qquad \qquad +  \lambda\varepsilon^{-1} \int \! \! \! \int_{Q_{\omega}} \phi^{1/2} |\nabla w| \xi \phi^{1/2} w dxdt  \biggl).  \nonumber
\end{split}\end{equation*}
Therefore, by  Young's inequality, we have that
\begin{align}
&  s\lambda^2 \delta^2\varepsilon^{-2}\int \! \! \! \int_{Q_{\omega}}\phi \xi^2|\nabla w|^2 dxdt \nonumber 
  \\&\qquad \leq  \frac{1}{4}\int \! \! \! \int_{Q} |L_2w|^2dxdt +   C\beta^{-2}s^3\lambda^4 (\delta^4 + \delta^2)	\varepsilon^{-2} \int \! \! \! \int_{Q_{\omega}}\phi^3|w|^2dxdt. \nonumber 
\end{align}
Thus, inequality  (\ref{eq31}) gives
\begin{equation}\begin{split} \label{eq35}
& ||L_1w||_{L^2(Q)}^2 + ||L_2w||_{L^2(Q_T)}^2 +  \beta^{-2}s^3\lambda^4\varepsilon^{-2}\int \! \! \! \int_{Q}\phi^3|w|^2dxdt   \\
& \qquad \qquad  \qquad \qquad \qquad \qquad \qquad \qquad +   s\lambda^2\varepsilon^{-2}\int \! \! \! \int_{Q}\phi |\nabla w|^2dxdt   \\
 &\qquad \leq C\biggl(  ||e^{s\alpha}g||_{L^2(Q)}^2+ \beta^{-2} s^3\lambda^4\varepsilon^{-2}\int \! \! \! \int_{Q_{\omega}} \phi^3|w|^2dxdt    \biggl). 
 \end{split}\end{equation}

Let us now we use the first two terms in left-hand side of (\ref{eq35}) in order to add the integrals of $|\Delta w|^2$ and $|w_t|^2$ to the left-hand side of  (\ref{eq35}). This is done using the expressions of $L_1w$ and $L_2w$. Indeed, from  (\ref{L1}) and (\ref{L2}), we have 

\begin{equation}\begin{split}\label{eq41}
& \int \! \! \! \int_{Q} s^{-1}\phi^{-1}|\pt w|^2 dxdt  + \varepsilon^{-2}\int \! \! \! \int_{Q} s^{-1}\phi^{-1}
\sum_{i,j= 1}^N \left | \pxi(M_{ij} \pxj w) \right | ^2 dxdt \\
&\qquad+  s^3\lambda^4\varepsilon^{-2}\int \! \! \! \int_{Q}\phi^3|w|^2dxdt   
+   s\lambda^2\varepsilon^{-2}\int \! \! \! \int_{Q}\phi |\nabla w|^2dxdt  \\
 & \leq C\biggl(  ||e^{s\alpha}g||_{L^2(Q)}^2+ s^3\lambda^4\varepsilon^{-2}\int \! \! \! \int_{Q_{\omega}} \phi^3|w|^2dxdt   \biggl).
\end{split}\end{equation}
Using the term in $|\pxi(M_{ij} \pxj w)|^2$ on the lef-hand side of (\ref{eq41}) and elliptic regularity, it is easy to show that 
\begin{eqnarray*}
s^{-1}\varepsilon^{-2}\int \! \! \! \int_{Q}\phi^{-1} \sum_{i,j= 1}^N \left | \pxij w  \right |^2dxdt
&\leq&C\biggl(  ||e^{s\alpha}g||_{L^2(Q)}^2+ s^3\lambda^4\varepsilon^{-2}\int \! \! \! \int_{Q_{\omega}} \phi^3|w|^2dxdt   \biggl).
\end{eqnarray*}
Estimate (\ref{eq41}) then gives 

\begin{eqnarray} \label{eq43}
&& \int \! \! \! \int_{Q} s^{-1}\phi^{-1}|\pt w|^2 dxdt  
+ s^{-1}\varepsilon^{-2}\int \! \! \! \int_{Q}\phi^{-1} \sum_{i,j= 1}^N \left | \pxij w \right |^2dxdt \nonumber \\
&&+  s^3\lambda^4\varepsilon^{-2}\int \! \! \! \int_{Q}\phi^3|w|^2dxdt   
+   s\lambda^2\varepsilon^{-2}\int \! \! \! \int_{Q}\phi |\nabla w|^2dxdt    \nonumber \\
 && \leq C\biggl(  ||e^{s\alpha}g||_{L^2(Q)}^2+ s^3\lambda^4\varepsilon^{-2}\int \! \! \! \int_{Q_{\omega}} \phi^3|w|^2dxdt   \biggl).
\end{eqnarray}
From (\ref{eq43}) and the fact that $w= e^{s\alpha}v$,  we finish the proof of Theorem \ref{carleman-calor}.
\end{proof}

\subsection{A Slightly changed Carleman inequality }
Our second sharp Carleman inequality is the following.
\begin{theorem}
For any $0< \epsilon \leq 1$, there exist $\lambda_0 = \lambda_0(\Omega, \omega)  \geq 1 $ and $s_0 = s_0(\Omega, \omega) \geq 1$ such that, for every $ \lambda \geq  \lambda_0$ and   $s \geq s_0 (T + T^2)$,  the solution $v$ of  \eqref{eq1-1} satisfies 
\begin{equation} \begin{split}  \label{ineq8}
&\int \! \! \! \int_{Q} e^{2s\alpha}|\pt v|^2 dxdt   + \varepsilon^{-2}\int \! \! \! \int_{Q}e^{2s\alpha} \sum_{i,j= 1}^N \left |\pxij v \right | ^2dxdt \\  &\qquad \qquad+ s^4\lambda^4\varepsilon^{-2}\int \! \! \! \int_{Q}\phi^4e^{2s\alpha}|v|^2dxdt   
+   s^2\lambda^2\varepsilon^{-2}\int \! \! \! \int_{Q}\phi^2e^{2s\alpha} |\nabla \rho|^2dxdt    \\
 & \leq C e^{\lambda||\psi||}( s\int \! \! \! \int_{Q}\phi e^{2s\alpha}|g|^2dxdt +  s^4\lambda^4\varepsilon^{-2}\int \! \! \! \int_{Q_{\omega}} \phi^4e^{2s\alpha}|v|^2dxdt),
\end{split}\end{equation}
with  $C>0$  depending only on $\Omega$, $\omega_0$, $\psi$ and $\beta$.
\end{theorem}

\begin{proof}
The starting point is the application of  the  Carleman inequality given in Theorem \ref{carleman-calor} to the equation \eqref{eq1-1}. Indeed, we have 

\begin{equation} \begin{split}\label{ineqB2}
&\varepsilon^2\int \! \! \! \int_{Q} s^{-1}\phi^{-1}e^{2s\alpha}|v_t|^2 dxdt  +\int \! \! \! \int_{Q}s^{-1}\phi^{-1}e^{2s\alpha} \sum_{i,j= 1}^N \left | \pxij v \right |^2dxdt 
\\  &\qquad \qquad + s^3\lambda^4\int \! \! \! \int_{Q}\phi^3e^{2s\alpha}|v|^2dxdt   
+   s\lambda^2\int \! \! \! \int_{Q}\phi e^{2s\alpha} |\nabla v|^2dxdt \\
 & \leq C(\int \! \! \! \int_{Q} e^{2s\alpha}|g|^2dxdt +  s^3\lambda^4\int \! \! \! \int_{Q_{\omega}} \phi^3e^{2s\alpha}|v|^2dxdt).
\end{split}\end{equation}

Next, we introduce the function $y(x,t) = v (x,t) (\phi^*(t))^{\frac{1}{2}}$. This new function satisfies
\begin{equation}\label{eqB4}
\left \{
\begin{array}{ll}
\varepsilon \pt y -\Div (M(x)\nabla  y) = -\varepsilon \frac{(T-2t)}{2}\phi^* y  +\varepsilon (\phi^*(t))^{\frac{1}{2}}g  &\qquad  \mbox{ in } Q, \\
y = 0 &\qquad  \mbox{ on } \Sigma.
\end{array}
\right. 
\end{equation}
Applying again the  Carleman inequality given by Theorem \ref{carleman-calor}, this time for $y$, we obtain, for $s$ large enough, that
\begin{equation} \begin{split} \label{ineq444}
&\int \! \! \! \int_{Q} s^{-1}\phi^{-1}e^{2s\alpha}|\pt y|^2 dxdt  + \varepsilon^{-2}\int \! \! \! \int_{Q}s^{-1}\phi^{-1}e^{2s\alpha} \sum_{i,j= 1}^N \left |  \pxij y \right|^2dxdt  \\  &\qquad \qquad + s^3\lambda^4\varepsilon^{-2}\int \! \! \! \int_{Q}\phi^3e^{2s\alpha}|y|^2dxdt   
+   s\lambda^2\varepsilon^{-2}\int \! \! \! \int_{Q}\phi e^{2s\alpha} |\nabla y|^2dxdt \\
 & \leq C( \int \! \! \! \int_{Q}\phi^* e^{2s\alpha}|g|^2dxdt  +   s^3\lambda^4\varepsilon^{-2}\int \! \! \! \int_{Q_{\omega}} \phi^3e^{2s\alpha}|y|^2dxdt).
\end{split}\end{equation}
From the definition of $y$, it is easy to show that
\begin{align}\label{sfdfdsf}
\int \! \! \! \int_{Q} s^{-1}\phi^{-1}e^{2s\alpha}| v_t (\phi^*)^{\frac{1}{2}}|^2 dxdt \leq
\int \! \! \! \int_{Q} s^{-1}\phi^{-1}e^{2s\alpha}|\pt y|^2 dxdt + \int \! \! \! \int_{Q} e^{2s\alpha}\phi |y|^2 dxdt. 
\end{align}
 Using \eqref{sfdfdsf}, inequality (\ref{ineq444}) becomes 
\begin{equation} \begin{split} \label{ineq6+}
&\int \! \! \! \int_{Q} s^{-1}\phi^{-1}\phi^*e^{2s\alpha}|v_t|^2 dxdt   + \varepsilon^{-2}\int \! \! \! \int_{Q}s^{-1}\phi^*\phi^{-1}e^{2s\alpha} \sum_{i,j= 1}^N \left | \pxij v \right |^2dxdt  \\  
&\qquad \qquad+ s^3\lambda^4\varepsilon^{-2}\int \! \! \! \int_{Q}\phi^3\phi^*e^{2s\alpha}|v|^2dxdt   
+   s\lambda^2\varepsilon^{-2}\int \! \! \! \int_{Q}\phi \phi^*e^{2s\alpha} |\nabla v|^2dxdt   \\
 & \leq C(\int \! \! \! \int_{Q}\phi^* e^{2s\alpha}|g|^2dxdt +  s^3\lambda^4\varepsilon^{-2}\int \! \! \! \int_{Q_{\omega}} \phi^3\phi^*e^{2s\alpha}|v|^2dxdt).
\end{split}\end{equation}
From Remark \ref{obs2}, the result follows.
\end{proof}

\section*{Acknowledgments} 
This paper has been partially established during the visit of M. Bendahmane to the Basque Center
for Applied Mathematics. The authors thank professor Enrique Zuazua for various fruitful  discussions about this work and Erich Foster for a careful reading of this paper.

\end{document}